\documentclass[12pt]{amsart}
\usepackage{graphicx}
\usepackage{amssymb}
\usepackage{epstopdf}
\theoremstyle{plain}

 \theoremstyle{definition}
 
\numberwithin{equation}{section}

\title{  DEGREE $k$ LINEAR RECURSIONS mod(p) AND NUMBER FIELDS}
\author{T. MacHenry and Kieh Wong}

\begin{document}




\maketitle

\begin{quote} ABSTRACT:  Linear recursions of degree $k$ are determined by evaluating the sequence of Generalized Fibonacci Polynomials, $\{F_{k,n}(t_1,...,t_k)\}$  (isobaric reflects of the complete symmetric polynomials)  at the integer vectors $(t_1,...,t_k)$. If $F_{k,n}(t_1,...,t_k) = f_n$, then 
 $$f_n - \sum_{j=1}^k t_j f_{n-j} = 0,$$ 
 and $\{f_n\}$ is a linear recursion of degree $k$.
 On the one hand, the periodic properties of such sequences modulo a prime $p$ are discussed, and are shown to be related to the prime structure of certain algebraic number fields; for example, the arithmetic properties of the period are shown to characterize ramification of primes in an extension field. On the other hand, the structure of the semilocal rings associated with the number field is shown to be completely determined by Schur-hook polynomials.\\
 key words:  Symmetric polynomials, Schur polynomials, linear recursions, number fields.
 \end{quote}

\vspace{1.0cm}

{\large{\emph{1.INTRODUCTION}}}

A sequence $\{f_n\}$ is a \textit{linear recursion of degree $k$}, denoted by $[t_1, ...,t_k]$,
if, given a (finite) sequence of integers $t_1, ...,t_k$, the following equation is satisfied for all $n\in \mathbb{\mathrm{Z}}$: 

$$f_n - \sum_{j=1}^k t_j f_{n-j} = 0.$$
In this paper we shall discuss the periodic nature of such sequences, and the periodic nature of such sequences modulo primes.  In particular,  we characterize those k-linear sequences which are periodic, and those which are periodic modulo a prime.  While we believe that these results are new and interesting,  it is the setting that they occur in, and the applications of these results that we are most interested in.  The setting in question is that of the ring of symmetric polynomials,  and the applications are to the theory of algebraic number fields on the one hand and to the theory of multiplicative arithmetic functions on the other. It is the first of these applications, the number fields,  that will be emphasized in this paper, while the second, the multiplicative arithmetic functions, will be discussed in more detail in a paper to follow shortly.

In [11], [12],[13],[14], the notion of isobaric polynomials was introduced, or rather reintroduced.  These are just the symmetric functions written in the Elementary Symmetric Polynomial (ESP) basis.  Historically, interest in the symmetric polynomials arose because of the famous relation between the roots of a (say, monic) polynomial and its coefficients.  If, for example,  we  take our monic polynomial to be
 $$X^k -t_1 X^{k-1} -...-t_k$$ 
 with roots  $$\lambda_1,...,\lambda_k$$
 then it is the classical result that the $t_j$ are, up to sign, the ESP's of the $\lambda_i$.  Regarding the $\lambda_i$ as indeterminants,  the $k-$degree symmetric   polynomials are those polynomials on the $\lambda_i's$,  which are invariant under the action of the symmetric group of degree $k$ acting on the generators.   All of this is, of course classical knowledge;  however,  when we rewrite the symmetric functions in the ESP basis, using the famous theorem that this is indeed possible,  we see that the form of the polynomials no longer emphasizes the symmetry of the generators,  but rather it is the partitions of the natural numbers which come to the fore.  After the mapping $$t_j = (-1)^{j+1}\mathcal{E}_j,$$ where $\mathcal{E}_j$ is the $j-th$ Elementary Symmetric Polynomial in $k$ variables, an \textit{isobaric } polynomial looks like this: 
 $$P_{k,n} = \sum_{\alpha\vdash n}C_\alpha t_1^{\alpha_1}...t_1^{\alpha_k}$$
 where  $ \mathbf{\alpha} = (\alpha_1,...,\alpha_k) $   is an integer vector with $\sum_{j=1}^k j\alpha_j = n$; that is,  $(1^{\alpha_1},...,k^{\alpha_k})$  is a partition of  $n$ into parts with   $\alpha_j \, j's $. We shall say that a symmetric polynomial written in this way,  emphasizing the partitions of the integers, has \textit{isobaric degree} $n$.  It can be thought of as a polynomial whose variables are Young diagrams (e.g. see [10]), or more accurately, the Young diagrams representing partitions of  $n$ into parts not larger than  $k$. Note that the coefficients  $C_\alpha$ are integers.
 
 Of special interest to us in this paper are the sequences of isobaric polynomials which form linear recursions, that is, sequences for which, given the variables $\mathbf{t} = (t_1,...,t_k)$,  we have for each $k$ a sequence of polynomials $\{P_{k,n}\}$ for which $$P_{k,n} = t_1 P_{k,n-1} + ... + t_kP_{k.n-k}.$$
 
 The mapping from the $\lambda-$basis to the ESP basis is a ring isomorphism,  that is,  we can speak of the (k-graded) ring of isobaric polynomials.  Letting $t_j=0$ for $j > k$ yields a projection of the ring onto the $k-th$ level of the grading.
  
 It is proved in [13] that such sequences form a free $k-$graded  $\mathbb{Z}-$module with a basis consisting of Schur-hook polynomials  [15] .
 In [13]  this was called the module of \textit{Weighted Isobaric Polynomials} or the WIP-module.  It can be thought of as a module of polynomials,  but is best considered as a module of sequences of polynomials.  In Section 3,  we suggest a way of looking at this module which is both intuitively transparent and algebraically very useful.  But first  we  discuss two especially important sequences in this setting.  They are the Generalized Fibonacci Sequence (GFP), and the Generalized Lucas Sequence (GLP) [12].  In the $\lambda-$basis, they are better known as the sequence of Complete Symmetric Polynomials,  and the sequence of Power Symmetric Polynomials.  
 
In the isobaric basis the GFP's are of the form, 
 
$$F_{k,n} = \sum_{\alpha\vdash n} \left(\begin{array}{ccc}  \;\;  |\alpha|  \\ \alpha_1 ... & \alpha_k\end{array}\right) t_1^{\alpha_1}...t_k^{\alpha_k},$$  

where $\alpha = (\alpha_1,...,\alpha_k), |\alpha| = \sum_{j=1,...,k} \alpha_j;$

 and the GLP's, of the form,
$$G_{k,n} = \sum_{\alpha\vdash n}\frac{n}{|\alpha|} \left(\begin{array}{ccc}  \;\;  |\alpha|  \\ \alpha_1 ... & \alpha_k\end{array}\right) t_1^{\alpha_1}...t_k^{\alpha_k}.$$  

In general, a Weighted Isobaric Polynomial (or WIP-polynomial) is given by the expression:

$$P_{\omega,k,n} = \sum_{\alpha\vdash n} \left(\begin{array}{ccc}  \;\;  |\alpha|  \\ \alpha_1, ... , \alpha_k\end{array}\right)\frac{\sum_j^k \alpha_j\omega_j}{|\alpha|} t_1^{\alpha_1}...t_k^{\alpha_k},$$  

where $\mathbf{\omega} = (\omega_1,...,\omega_k)$ is a \textit{weight} vector. 
Each weight vector determines a k-linear recursion [12]. The weight vectors for the generalized Fibonacci sequence and the generalized Lukas sequence are given, respectively,  by  $\mathbf{\omega} = (1,1,...,1,...)$,  and by  $\mathbf{\omega} = (1,2,...,n,...)$.  It is straightforward to check that $P_{{\omega},k,n}$ is a k-linear recursive sequence of isobaric polynomials for each $\omega$ and $k$.
And,  that we can add two weighted sequences in the WIP-module by adding their weight vectors, thus realizing the abelian group structure of the module and emphasizing that the preferred basic element of the module is a sequence.

For each $k$ we call any weighted sequence of isobaric polynomials a \textit{generic} $k-$linear recursive sequence, allowing the application of an evaluation map to the indeterminates.  The WIP-module, and, indeed,  the entire ring of $k-$isobaric polynomials $P_{k,n}$ is determined implicitly by the $k-$degree monic polynomial 
$$\mathcal{C}(X) = X^k-t_1X^{k-1}-...-t_k = \sum_{j=1}^k  t_j X^{k-j}$$

by virtue of the two fundamental theorems of symmetric functions alluded to above,  and the change of basis. Therefore, we call this polynomial 
$$\mathcal{C}(X) =  \sum_{j=1}^k  t_j X^{k-j},$$ 
the \textit{Core Polynomial}.  Thus given the core polynomial,  the whole isobaric structure falls into place. But the core polynomial itself is uniquely given once we have assigned the generic variables $t_1,...,t_k$,  so we find it convenient to use the notation $[t_1,...,t_k]$ to denote the core
polynomial.  Since we can take $k$ to be arbitrarily large,  it is also convenient to give power series the honourary status of core polynomial (with some adjustment necessary to the bracket notation).

These remarks will be more effective when we look not just at the generic core,  but also consider evaluation maps on the $\mathbf{t}-$vectors, that is when we look at polynomials of degree $k$ with numerical coefficients. (In section 3, the unity of these ideas will become especially transparent).

Suppose we choose to evaluate the indeterminates $\mathbf{t}$ in the ring of integers, then each sequence  $\{P_{\omega,k,n}(\mathbf{t})\}$ gives a numerical  $k$-degree linear recursion;  and since, in particular, $ \mathbf{t}$ is given, the core is uniquely determined.  Moreover,  every k-degree linear recursion can be realized in this way.  The contents of Section 3 will suggest that choosing to use the GFP as our  generic sequence,  has a great deal of merit. This sequence contains the polynomials

\begin{enumerate}
  \item $F_{k,0} = 1$
  \item $F_{k,1} = t_1$
  \item $F_{k,2} = t_1^2 + t_2$
  \item $F_{k,3} = t_1^3 + 2t_1t_2 + t_3$
  \item $F_{k,4} = t_1^3 + 3t_1^2t_2 + t_2^2 + 2t_1t_3 + t_4 $
  \item etc. 
\end{enumerate}

of isobaric degrees $0, 1,2,3,4$.

If we let $k=2$, that is,  use the projection $t_j = 0 $ for $j>2$, and let $[t_1,t_2] = [1,1] $, we find that the sequence  $\{F_{2,n}\}$ is just the Fibonacci sequence.  A similar exercise for the GLP's yields the Lucas sequence.  In either case,  the core polynomial is $X^2 -X - 1$.  However, once a core polynomial is chosen, given the $k$  and the generic linear recursion, all is determined. So in particular,  if we choose a generic k-linear recursion, then for each evaluation of the 
$\mathbf{t}-$vector, exactly one core polynomial is selected.  In this way, we get a one-to-one relation between k-cores and all numerical linear recursions. 

So the first two questions we ask are:

(1) Which  numerical linear recursions are periodic? 

and 

(2) Which are periodic mod(p)? ( A sequence $\{f_n\}$ is periodic if there is a positive integer $c$ such that for all $n$,  $f_{n+c} = f_n$ .)  

At the end of this paper,  we include a Maple Algorithm,  due to Professor Mike Zabrocki of York University, for computing the period of any $k$-order linear recursion modulo a prime $p$. 

{\large {\emph{2.PERIODIC LINEAR RECURSIONS}}}

We now answer the two questions asked in Section 1, reminding the reader that the generic recursion that we are using is the GFP sequence.  While these results are not  difficult to prove, it seems that they are not stated in the literature. 
 
 $\mathbf{THEOREM  \,2.1}$

A linear recursion is periodic if and only if every root of the core polynomial is a primitive complex root of unity.  In particular, if the core polynomial is the cyclotomic polynomial $CP(n)$ of degree  $\phi (n)$, where $\phi$ is the Euler totient function, then its associated linear recursion is periodic with period $n$ \cite{MT}.  (It is interesting to compare this theorem with the Lech-Mahler Theorem \cite{Cassels}.) $\square$
 
The proof will be discussed in Section 3.

Denote the period of a linear recursion, either mod(p) or mod(1),  by  $c_p[t_1,...,t_k]$ where $p$ is either a rational prime or $p=1$, and the $t_j$ are the coefficients of the core polynomial.

$ \mathbf{THEOREM \,2.2}$

Every linear recursion is periodic modulo  $p$ for every rational \textit{prime} $p$. The period  $c_p[\mathbf{t}] \leqslant p^k$. $\square$

This follows from simple combinatorial arguments, essentially the pigeonhole principle.  We improve this bound to a best bound, $p^k - 1$,  in Section 4. We observe that if a sequence is a periodic linear recursion, then
$F_{k,c_p} = F_{c_p} = 1, F_{c_{p}-1} = ...= F_{c_p -k+1} = 0,\, F_{c_p+1} = t_1.$


While there is nothing deep about the proofs of these two theorems,  it is of some interest that they occur within the confines of the ring of symmetric functions and become obvious when this particular basis is chosen.  However, the particular techniques for studying recursions and periodic recursions reveal an even deeper connection with symmetric functions framed in the language of isobaric polynomials,  which in turn points to a strong connection with combinatorial algebra.  In the next section we discuss some not so well-known "well-known" results, and add some new information which we believe not to be well-known. We now discuss our most important tool;  namely, the companion matrix of the core polynomial and a rather remarkable structure induced by it.

3.$\large{\emph{THE COMPANION MATRIX OF THE CORE POLYNOMIAL}}$

With each core polynomial, we associate its rational canonical matrix, the so-called \textit{companion matrix}.  We first consider the companion matrix for the generic core polynomial of degree $k$.

\begin{displaymath}
\mathbf{A} = 
\left(\begin{array}{cccc}
0&1&...&0\\
0&0&...&0\\
0&0&...&1\\
t_k&t_{k-1}&...&t_1
\end{array} \right) 
\end{displaymath}

Since $det{\bf A}= (-1)^{k+1} t_k$, $det{\bf A^n} = (-1)^{n(k+1)} t_k^n,  $ ${\bf A}$ is singular iff $t_k =0$. But if $t_k = 0$, the core polynomial is reducible;  so we assume  $\mathbf{A}$ to be non-singular.  
Thus $\mathbf{A}$ is invertible and generates a cyclic group (finite, if the coefficients of the core polynomial satisfy the conditions of Theorem 2.1, otherwise, infinite).   The inverse of $\mathbf{A}$  is 


\begin{displaymath}
\mathbf{A^{-1}} = 
\left(\begin{array}{ccccc}
-t_{k-1}t_k^{-1}&-t_{k-2}t_k^{-1}&...&t_1t_k^{-1}&t_k^{-1}\\
1&0&...&0&0\\
0&1&...&0&0\\
...&...&...&...&...\\
0&0&...&1&0
\end{array} \right).
\end{displaymath}

We record the orbit of the k-th row vector of $\mathbf{A}$ under the action of $\mathbf{A}$, below $\mathbf{A}$, and the orbit of the first row of $\mathbf{A}$ under the action of $\mathbf{A^{-1}}$ on the first row of $\mathbf{A}$ is recorded above $\mathbf{A}$, and consider the $ \infty \times k$ matrix whose row vectors are the elements of the doubly infinite orbit of  $\mathbf{A}$ acting on any one of them. 
For $k=3,$ $\mathbf{A^\infty}$ looks like this (we explain the symbols for the elements below):
\begin{displaymath}
\mathbf{\mathbf{A^\infty}} =
\left(\begin{array}{ccc}
...&...&...\\
S_{(-n,1^2)} & -S_{(-n,1)} & S_{(-n)} \\
...&...&...\\
S_{(-3,1^2)} & -S_{(-3,1)} & S_{(-3)}\\
1&0&0\\
0&1&0\\
0&0&1\\
t_3&t_2&t_1\\
... & ... & ...\\
S_{(n-2,1^2)} & -S_{(n-2,1)} & S_{(n-2)}\\
S_{(n-1,1^2)} & -S_{(n-1,1)} & S_{(n-1)}\\
S_{(n,1^2)} & -S_{(n,1)} & S_{(n)}\\
... & ... & ...

\end{array}\right)_{\infty \times 3}
\end{displaymath}

\noindent This matrix has a number of important features which we summarize in

 $\mathbf{THEOREM \, 3.1}$(cf.\cite{HM}, \cite{Lascoux}, \cite{Lascoux2}, \cite{CV1}, \cite{CV2})

(3.11) The row vectors  consist of the orbit of any row with $\mathbf{A}$ acting as a transformation matrix (on the right, say), and the components of the row vectors are just isobaric reflects of Schur-hook polynomials.

(3.12) The set  of $k\times k$ contiguous row vectors of $\mathbf{A^\infty}$,  with the entry in the lower right hand corner being $\mathbf{S_{(n)}}$, yields a (faithful) matrix representation of the cyclic group generated by $\mathbf{A}$:

  \begin{displaymath}
\mathbf{\mathbf{A^n}} = 
\left(\begin{array}{ccccc}
(-1)^{k-1}S_{(n-k+1,1^{k-1})}&...&(-1)^{k-j}S_{(n-k+1,1^{k-j})}&...&S_{(n-k+1)}\\
...&...&...&...&...\\
(-1)^{k-1}S_{(n,1^{k-1})}&...&(-1)^{k-j}S_{(n,1^{k-j})}&...&S_{(n).}
\end{array} \right)
\end{displaymath}

Or, more succinctly, we have  $${\bf A^n} = [(-1)^{k-j} S_{(i,1^{k-j})}]_{k\times k},$$
where the entries are  isobaric Schur-hook reflects whose Young diagrams have arm length $i$ and leg length $k-j$ in the case of positive $n$.  

(3.13) The elements in each row of $\mathbf{A^\infty}$ are the coefficients of a representation of the powers (positive and negative) of any of the roots of the core polynomial---denoted by $\lambda^n$--- in terms of a basis consisting of the first $k-1$ powers of  $\lambda$:  $$\lambda^n = \sum_{j=0}^{k-1}(-1)^{k-j}S_{(n,1^{k-j}) }\lambda^{j}$$  for $n \in\mathbb{ Z}$,  where $\lambda$ is a root of the core polynomial  (and, as remarked above,  the coefficients are Schur-hook reflects whose Young diagrams have  arm length $n$ and leg length $k-j$ when  $n$ is positive). 

(3.14) Each column of $\mathbf{A^\infty}$ is a $\textbf{t}$-linear recursion of Schur-hook polynomials. In particular, the right hand column is just the (doubly infinite) sequence of Generalized Fibonacci Polynomials, $\mathbf{F_{k,n}}$.

(3.15) $tr(\mathbf{A^n}) = \mathbf{G_{k,n}}(\textbf{t}$) for $n \in \mathbb{Z}$, where $\mathbf{G_{k,n}}$ is just the sequence of Generalized Lucas Polynomials, which is also a \textbf{t}-linear  recursion.

(3.16)  It is appropriate also to call the negatively indexed entries in the matrix Schur-hook polynomials.  It is of interest that they can be represented as quotients of two positively-indexed Schur polynomials,  which in general are not hook polynomials:  For example. $k=3$ is a typical case: 

\begin{center}
 
  (1) $S_{(-n)}$ = $\frac{S_{{(n-3)^2})}}{t_{3}^{(n-2)}};$\\ 
  (2) $S_{(-n,1)}$ = $-\frac{S_{{(n-2,1)})}}{t_{3}^{(n-2)}};$ \\
  (3) $S_{(-n,1^2)}$ = $\frac{S_{{(n-2)^2})}}{t_{3}^{(n-2)}}$ \\

\end{center}

\textbf{REMARK}  We note that the existence of the matrix $\mathbf{A^\infty}$ extends the sequences of Schur-hook polynomials, in particular, the GFP, as well as the GLP, in the negative direction.  It would be interesting to have a combinatorial interpretation of these negatively indexed symmetric functions.
One might compare this result with the theorem in \cite{TM2}, which gives rational convolution roots to all of the elements in the WIP-module \cite{MT2}, i.e., to all of the sequences of symmetric functions in the free $\mathbf{Z}$-module generated by the Schur-hook polynomials. 

Proofs (3.11-3.15).

(3.11) The orbit structure is a consequence of the construction of the matrix. Operation of the companion matrix on a $k$-vector of integers generates a linear recursion with respect to the vector $\mathbf{t}$. In fact,  the Schur-hook sequences sequences claimed in the theorem  \cite{MT}.

(3.12) follows from the arguments in (3.11).

(3.13)  follows from the Hamilton-Cayley Theorem.  A simple induction shows that these coefficients are just the stated Schur-hook functions of the theorem.

(3.14) This is discussed in (3.11). \\

(3.15) The traces of the $k \times k$ -blocks are the sums of all of the Schur-hook (reflects) whose Young diagrams partition the same n;  but such sums of Schur-hooks are well known to be GLP of isobaric degree $n$  \cite{MT}. $\square$

(3.16)   This can be easily proved using the recursion properties of the sequences; however, since this result will not be used in the paper,  we omit the proof.

The infinite companion matrix is a remarkable summary of all of the features connected with linear recursions (as enumerated in  Theorem 3.1):  It contains representations of the roots of the core polynomial as row vectors;  the right-hand column consists of GFP's, i.e., the generic k-th order linear recursions;  it displays the role of Schur-hook functions as both constituents of sequences of k-th order linear recursions---one of which is the GFP sequence--- and as  coefficients for a representation of the powers of the roots of the core polynomial. It contains a matrix representation of the free abelian group generated by the companion matrix, in particular,  a matrix representation of the free abelian group generated by any of the roots of the core.  It also contains, as traces, the GLP's.  Recalling that the GFP's and the GLP's are respectively, the isobaric versions of the complete symmetric polynomials and the power symmetric polynomials. With this we have shown a connection between the theory of linear recursion and an important submodule of the algebra of symmetric polynomials, the WIP-module. Moreover, we have introduced an extension of the symmetric polynomials to negatively indexed symmetric functions which are related to the reciprocals of powers of the roots of the core polynomial.  Thus we have a striking summary of the connection between the theory of equations and the theory of linear recursions within the ring of symmetric polynomials.  We note that while many of these properties of the extended companion matrix are known to A. Lascoux and his students [4],[5],[7],[8], the role of the GFP's and the GLP's, as well as the form of the negative entries, may not be so well-known. This matrix will be a useful and important tool in what follows.

$\mathbf{COROLLARY} \,3.2$ 

Given the k-th order linear recursion determined by the core $[t_1,\cdots,t_k]$,  with the companion matrix  $\mathbf{A}$, and denoting the cyclic group generated by $\mathbf{A}$ as $\mathbf{H}$,  we have that $\mathbf{H}$ is a finite cyclic group exactly when the linear recursion
is periodic, the order of the cyclic group $\mathbf{H}$ 
being the period of the recursion. Moreover, if the core polynomial is irreducible over the rationals, then every root of the core polynomial generates a finite cyclic group whose order is also the period of the linear recursion. 

Proof.   The proof follows immediately from Theorem 3.1. $\square$

Corollary 3.2 explains why Theorem 2.1 is true, for clearly, the only irreducible core polynomials having all of its roots periodic are those whose roots are roots of unity.

Applying the facts learned above about the companion matrix,  we now consider the periodic behaviour of linear recursions modulo a prime $p$. 
\vspace{0.50cm}

4.{\large{\emph{p-PERIODICITY AND THE COMPANION MATRIX}}}

Since the vector  $\mathbf{t} = [t_1,...,t_k]$ determines both the core polynomial and its associated linear recursion uniquely, we write  $[t_1,...,t_k]$ to denote either of these structures when the context is clear, and we shall extend the usage to the notation $[t_1,...,t_k]_p$ for a linear recursion $[t_1,...,t_k]$ modulo the prime $p$. As in Section 2, $c_p[t_1,...,t_k]$ denotes the period of $[t_1,...,t_k]_p$, and  $\mathbf{A_p}$ denotes the companion matrix with entries modulo $p$.  For any matrix $\mathbf{M}$, $tr\mathbf{M}$ denotes the trace of the $\mathbf{M}$.

THEOREM 4.1

(4.11)  $c_p[\mathbf{t} ] = c_p[t_1,...,t_k] \le p^k - 1. $
 
(4.12)  The (cyclic) group generated by $\mathbf{A_p}$ has order $c_p[\textbf{t}].$

(4.13) The columns of $\mathbf{A_p^\infty}$ have period $c_p[\textbf{t} ]$ .

(4.14) $\lambda^{c_p} [\textbf{t}] =_p 1$, where  $\lambda$ is a root of $\mathcal{C}{[\textbf{t}} ]$, and 
$c_p[\textbf{t]}$ is the least positive integer for which this is true;  i.e., $c_p[\textbf{t} ]$ is the $p$-order of $\lambda$.

(4.15) $tr\mathbf{A_p^n}$ is linearly recursive with period $c_p[\textbf{t}]$.

Proof (4.11-4.15).   (4.11) will be proved in  the next section.
Clearly $\mathbf{A_p}$ generates a cyclic group of order dividing $c_p$; on the other hand, since each of the columns of $\mathbf{A^\infty}$ is a linear recursion, they too must have a period $c_p$.
(4.14) is a direct consequence of (3.13),(3.14) and (4.13).
(4.15) is a consequence of (3.15). $\square$

\textbf{REMARK}:  As pointed out above, Corollary 3.2 accounts for the truth of Theorem 2.1.  The core polynomials for the primitive $n-th$ roots of unity are the cyclotomic polynomials of degree $\phi(n)$,  whose roots have the obvious geometric period of $n$; that is,  $c_p[\textbf{t}]$ = $n$,  where $\textbf{t}$ is the appropriate vector of coefficients of the cyclotomic polynomial of degree $\phi(n)$.  This also affords a geometric interpretation of periodicity for the roots of the core polynomial in the plane of complex numbers with coordinates taken mod$(p)$,  which is analogous to the cyclotomic periodicity.
\vspace{1cm}

5.{ \large{\emph{THE NUMBER FIELD $\mathbf{R}[\mathbf{t}]$ \textsc{AND THE SEMILOCAL RING }$\mathbf{R_p}[\mathbf{t}]$}}}

$\mathbf{PROPOSITION} \,5.1$
If the core polynomial  [\textbf{t}] is reducible $mod(p)$, and if  $p$ does not divide $c_p$, then $c_p[\textbf{t}]$ is the least common multiple of  the p-periods of its irreducible factors.  Otherwise, $c_p$ will be a proper multiple of the least common multiple of the  p-periods of its irreducible factors $\square$

\textbf{REMARK} Note that if the core polynomial is reducible,  it is reducible $mod(p)$.

\textbf{REMARK} Clearly,  there should be a better theorem which specifies the exact  "proper multiple" in the second sentence of Proposition 5.1.  That better theorem is best stated and proved in the context of number fields.  This is done in Section 6,  Theorem 6.8 and Corollary 6.9,  where this "proper multiple" in the statement of 5.1 is supplied.  Hence, we  defer the proof of Proposition 5.1 to Section 6.

Proposition 5.1 suggests that in many cases we may as well consider only irreducible cores. But in that case,  we can also consider the number field $\mathcal{F}$ = $ \mathbb{Q}(\lambda) $ = $\mathbb{Q}[X]/id<\mathcal{C}(X)>$.
Let us denote the ring of integers (the maximal order) in this field by $\mathbf{R}[\textbf{t}]$ and we write
 $\mathbf{R} \otimes \mathbb{Z}_p = \mathbf{R_p}$.  We can write the elements of the field $\mathcal{F}$ either as a module over the basis $\{1, \lambda,...,\lambda^{k-1}\}$, or uniquely as k-tuples $(m_o,...,m_{k-1})$ with entries from $ \mathbb{Q}$ with multiplication determined by the minimal polynomial of the field, or, as a result of the Hamilton-Cayley Theorem, as a module with the basis $\{\mathbf{I}, \mathbf{A},...,\mathbf{A^{k-1}}\}$. This gives a matrix representation of the elements in the field.  Call it the \textit{standard\textit{}} representation. We also have the same three options in $\mathbf{R_p}$ using these bases modulo $(p)$. Theorem (3.13) 
can be regarded as giving a representation of the powers of $\lambda$ in   $\mathcal{F}$, as polynomials in the integral $\lambda$ -basis where the coefficients are
Schur-hook polynomials evaluated at $[\textbf{t}]$. Note that we have an induced \textit{standard} matrix representation in the ring $\mathbf{R_p}$.

One of the concerns of the theory of algebraic number fields is the relation between primes in the extension field $\mathcal{F}$ and the rational primes in $\mathbb{Z}$ that they sit over.  If we let $p$ be a rational prime generating the prime ideal $\mathbf{p}$ in  $\mathbb{Q}$, and let $\mathcal{P}$ be the ideal in $\mathbf{R}$ extending $\mathbf{p}$, then $\mathcal{P} = \mathcal{P}_{1}^{\epsilon_{1}}...\mathcal{P}_{s}^{\epsilon_{s}}$  is the prime decomposition of $\mathcal{P}$ in the Dedekind ring $\mathbf{R}$. If $f_j$ is the relative degree of the prime ideal $\mathcal{P}_{j}$, i.e., the degree of its minimal polynomial, then either $s = 1$ and  $\epsilon_1 = 1$, in which case $\mathcal{P}$ is a prime ideal, and $p$ is \textit{inert};  or, $s >1$ but $\epsilon_j = 1$ for all $j's$, in which case  $\mathcal{P}$ is the product of distinct prime ideals, and  $p$ \textit{splits};  or,  some $\epsilon_j > 1$ and $p$ \textit{ramifies}. These properties are reflected in the semilocal ring $\mathbf{R}_p$. Moreover,  there is a relation between the  phenomenon of periodicity of the linear recursion associated with the core polynomial and properties of the primes in the extensions of the core localized at $p$.  This will be discussed in the following sections. It is well known that for each irreducible core polynomial 
only a finite number of primes ramify;  when they do,  they divide the discriminant of the field. With few exceptions, the converse is also true, and those exceptions will not occur in our discussion \cite{J};  hence, for the purposes of this paper, $p$ ramifies if and only if $p|\Delta$,  where $\Delta$ is the discriminant of $\mathcal{F}$.  We shall want to use the following well-known fact.
 $\mathbf{PROPOSITION}\, 5.2$

$\Delta$ = $(-1)^{k(k-1)/2} \mathbf{ N(\mathcal{C}(X))}det\mathcal{C^{\prime}}(\mathbf{t})$.\;\;\quad\quad $ \square$

$\mathbf{\mathcal{C}'}(\mathbf{t})$ is the derivative of the core polynomial, that is, the \textit{different}.

\noindent Noting that $ \mathbf{\mathcal{C'}}$  can be regarded as an element of $\mathbf{R_p}$, and, denoting  $ \mathbf{\mathcal{C'}}$ by $\mathbf{D}$, we have

$\mathbf{COROLLARY}\, 5.3$

$ \mathbf{D_p} $ generates an ideal in $\mathbf{R_p}$ (the discriminant ideal) if and only if $p|\Delta$, that is, if and only if $p$ ramifies in $\mathbf{R}$.

Proof.  $p$ divides the discriminant of the core polynomial  modulo $p$ if and only if $p$ ramifies,  which occurs if and only if the different vanishes modulo $p$ at a root of the core polynomial, and this happens if and only if the different generates an ideal in the semilocal ring $\mathbf{R_p}$ (the alternative being that the different is a unit in $\mathbf{R_p}$). $\square$

In keeping with the notation $\mathbf{A^\infty}[\mathbf{t}]$, we let $\mathbf{M^\infty}[\mathbf{t}]$ be the $\mathbf{H_p}$ -orbit of any row vector in the matrix $\mathbf{M}$. Since, by construction, the columns of a standard matrix are $\mathbf{t}$-linear recursions, the following proposition can be proved by induction.

\newpage$\mathbf{PROPOSITION}\, 5.4$

The right-hand column of $\mathbf{D^\infty}[\mathbf{t}]$ is the sequence of GLP's; that is,  the right hand column of $\mathbf{D^\infty}[\mathbf{t}]$ is a list of the traces of the matrices representing $\mathbf{A^n}$, thus, the right hand column of the matrix consists of the terms of the GLP-sequence  (cf. (3.15). 

\textbf{PROOF} If the core polynomial is $X^k - \sum_{j=1}^k t_j X^{k-j}$,  then the different  \textbf{D} is the polynomial $\textbf{D} = kX^{k-1} -  \sum_{j=1}^{k-1} t_j (k-j)X^{k-j-1}$.  This is represented by the vector $(-t_{k-1},...,-(k-1)t_1, k)$ in $\textbf{R}$, and the orbit of this vector under the action of the companion matrix of the core polynomial gives the standard matrix representation.  Since acting on a vector in  $\textbf{R}$ by $\textbf{A}$ automatically generates linearly recursive columns determined by $[t_1,...,t_k]$ it is necessary 
only to notice that the element in the  upper right hand corner of the matrix representing  $\textbf{D}$ is  k.  Induction does the rest.$\square$

\textbf{REMARK} There is an interesting connection between the GFP-sequence and the GLP-sequence; namely, they are related by partial differentiation.  Precisely,  $\frac{\partial}{\partial t_j} G_n = nF_{n-j}, j=1,...,k$ (cf. this with D.H. Lehmer's notion of companion sequences [9] ). 

\textbf{REMARK} It follows from Proposition 5.4 that the period of the different $ \mathbf{D_p} $ is the same as   $c_p(\mathbf{R_p})$ if p does not ramify.  If $p$ splits (recall that in this paper 'splits' mean factors but does not ramify) then $\mathbf{D_p}$ is in the group of units $\mathbf{G_p}$, and is a coset of $\mathbf{A_p}$, possibly identical with $\mathbf{A_p}$.  If $p$ ramifies,  then $\mathbf{D}$ is a maximal ideal in $\mathbf{R_p}$.  Theorem 5.4 and the remark above give the rather pretty set of connections among the GFP-sequence,  the core polynomial, the derivative of the core, and the GLP-sequences: GFP determines the core, the derivative of the core yields GLP, the derivative (any first partial) of GLP yields the GFP. 

6. {\large{\emph{STRUCTURE OF THE SEMILOCAL RING $\mathbf{R_p} = \mathbf{R\bigotimes\mathbb{Z}_p}$}}}

$\mathbf{R_p}$ is a finite, commutative ring;  it is, therefore, a semilocal ring.  The structure of semilocal rings is well-known (e.g., [6], [15]).  We restate the structure theorem here (Theorem 6.4) for 
easy reference.  $\mathbf{R_p}$ also has an orbit structure under the action of the group generated by $\mathbf{A_p}$, which, while not mysterious, is not readily found in the literature, and plays an integral role in our results. We shall first discuss this orbit structure and then exploit the semilocal nature of $\mathbf{R_p}$.

If $t_k \neq 0\, mod(p)$, then $\mathbf{A_p}$ is non-singular, and, hence, is a unit in $\mathbf{R_p}$. The units in  $\mathbf{R_p}$ are exactly those elements with norms different from $0$, that is, having a standard matrix with non-zero determinant. An element with zero norm, then, either is zero or belongs to a proper  ideal.  Denote the group of units of $\mathbf{R_p}$ by $\mathbf{G_p}$ and its subgroup $gp<\mathbf{A_p}>$ by $\mathbf{H_p}$, the \textit{period subgroup}. Then  $\mathbf{R_p}$ is a $\mathbb{Z}_p (\mathbf{H_p})$-module, or more conveniently,  a right $\mathbf{H_p}$-module.  Clearly, $\mathbf{R_p}$ is the disjoint union of its orbits under the action of $\mathbf{A_p}$. A number of observations follow from these facts. It will be useful to list them for future reference:

(1)  The orbit of zero is a singleton.\\
(2)  An ideal consists of the disjoint union of orbits, each of which has orbit length dividing $c_p[\mathbf{t}]$. (Clearly, two orbits are either disjoint or identical, up to cyclic permutation.)\\
(3) Two distinct orbits in the same maximal ideal differ from one another by a coset representative of  $\mathbf{H_p}$, i.e., if $O_1$ and $O_2$ are distinct orbits in the maximal ideal  $I$, then there is a coset representative $g$ of $\mathbf{H_p}$ in $\mathbf{G_p}$ such that $O_1g = O_2$. (Of course,  a coset representative may belong to the stabilizer of $\mathbf{H_p}$. $O_1$ and $O_2 $ need not be bijective.)\\
(4) The orbits of $\mathbf{G_p}$ are the cosets of $\mathbf{H_p}$.\\
(5) The columns of an orbit are $\mathbf{t}$-linearly recursive,  with a period dividing $c_p[\mathbf{t}]$.\\ 
(6) The (standard) matrix representation of $\mathbf{R_p}$ is implicit in the orbit structure of $\mathbf
{R_p}$.  

If $\mathbf{m} \in \mathbf{R_p}$, and if $m_{i,j}$ is the $(i,j)-th$ component of the standard matrix representation $\mathbf{M}_p$ of  $\mathbf{m}$, the row vectors $\mathbf{m_i}$ of $\mathbf{M_p}$ are just the elements of the $\mathbf{A_p}$-orbit  of $\mathbf{m}$. 
A pseudo-Hasse diagram that illustrates the construction of a typical finite ring $\mathbf{R_p}$  is given below: 
\newpage

\begin{figure}[h]
\centering
\includegraphics[width=1.5in]{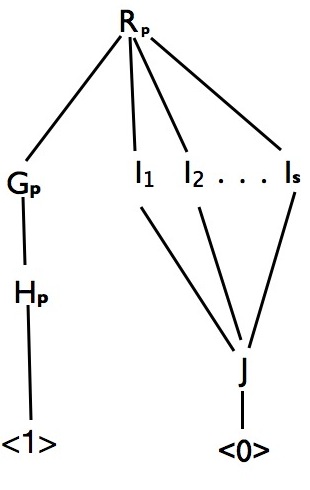} 
\caption{{\bf lattice diagram of semi-local ring}}
\label{semi-local ring}
\end{figure}
where  $\mathbf{G_p}$ is the group of units in the finite ring $\mathbf{R_p}$, $\mathbf{H_p }= gp<\mathbf{\mathbf{A_p}}>$, $|\mathbf{H_p}| = c_p$ = the period of the associate $k-$linear recursion $F_{k,n}(t_1,...,t_k) mod(p)$,  the $\mathbf{I_j}$ are the maximal ideals
 in this ring,  and $\mathbf{J}$ is the radical.  The "pseudo"  in pseudo-Hasse refers to the fact that we show the lattice structure of the group of units in this ring in the same diagram.  Here is an example:

\begin{figure}[h] 
   \centering
   \includegraphics[width=3in]{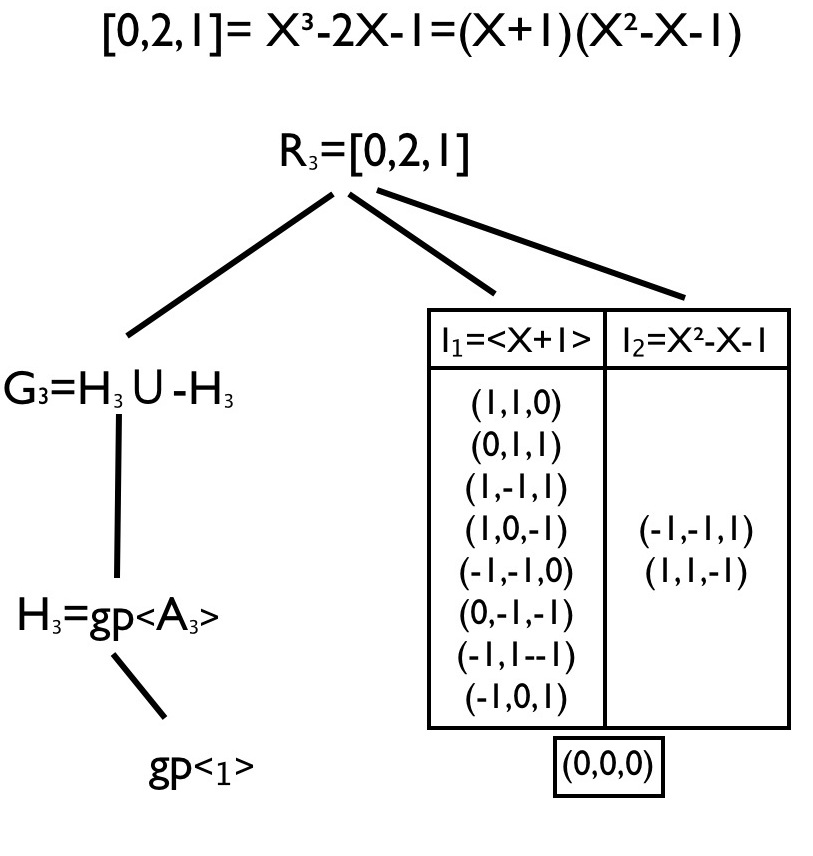} 
   \caption{{\bf lattice diagram of [0,2,1]}}
\end{figure}

$\mathbf{e_1} = (-1,1,-1), \mathbf{e_2} = (-1,-1,1), |\mathbf{R_3}| = 27, |\mathbf{G_3}| = 16, |\mathbf{H_3}| = 8, |\mathbf{G_3:H_3}| = 2.$

\begin{center}
$ \mathbf{A_3 }= \left(\begin{array}{ccc}  0 & 1  & 0  \\ 0 & 0  & 1  \\ 1 & -1  & 0 \end{array}\right)$
\end{center}

We record the following well-known fact:

$\mathbf{PROPOSITION } \,6.1$

There is a one-to-one correspondence between  maximal ideals of $\mathbf{R_p}$ and  irreducible factors of $\mathcal{C }(X) mod(p) $.  $\square$

$\mathbf{PROPOSITION }\,6.2$ (Traces)

Let $\mathbf{m} = (m_0,...,m_{k-1}) \in \mathbf{R_p}$.  $tr(\mathbf{m}) = m_0 \mathbf{G_{k,0}} + ... +m_{k-1} \mathbf{G_{k,k-1}} $, where $\mathbf{\{G_{k,n}\}}$
is the sequence of Generalized Lucas Polynomials, i.e, the isobaric reflect of the complete symmetric polynomials. 

Proof. Express $\mathbf{m}$ as $(m_0,m_1,...,m_{k-1})$ and note that the rows of $\mathbf{M}$ are vectors $\mathbf{mA^i}$. Writing $\mathbf{A_j^i}$  for the $j-th$ column of $\mathbf{A^i}$, we have that the trace of $\mathbf{m}$ is
$$(\mathbf{mA^0}) \mathbf{A_1^0} + (\mathbf{m }\mathbf{A^1}) \mathbf{A_2^1} + ... +(\mathbf{m}\mathbf{A^{k-1}}) \mathbf{A_k^{k-1}}.$$
But a suitable rearrangement of this sum is just
$$m_0 tr\mathbf{A^0} + m_1 tr\mathbf{A^1} + ...+m_{k-1} tr\mathbf{A^{k-1}} .$$
which, by Theorem 3.15, yields Theorem 6.2. $\square$

Also note that since each component of a vector in an orbit is in exactly one trace computation, the sum of the components of vectors in an orbit is equal to the sum of the traces of the vectors in the orbit.  That is,

$\mathbf{PROPOSITION} \,  6.3$

(6.31) The sum of the elements of the $\mathbf{A_p}$ -orbit  of the vector $\mathbf{m}$ is the sum of the traces of the row vectors, $\mathbf{m_i}$, i.e.,
$$ \sum_{i,j} m_{i,j}= \sum_i tr \mathbf{m_i}.$$

(6.32)  If $\mathbf{m }\in \mathbf{R_p}$ , i.e, if  $det \mathbf{m}$ = $0$, then $$\sum_{orbits \, of\, I} \sum_i tr (\mathbf{m_i}) = 0. \qquad\qquad\square$$ 

\newpage $\mathbf{THEOREM} \,6.4$ (e.g.,\cite{McD}VI.2, \cite{J}). 

$\mathbf{R_p}$ is a semilocal ring. In particular, letting $\mathbf{J}(\mathbf{R_p}) = Rad(\mathbf{R_p}) = \mathbf{I_1} \cap ... \cap \mathbf{I_s }= \mathbf{I_1} ...\mathbf{I_s}$,  $\mathbf{I_1}, ...,\mathbf{I_s}$  a complete set of maximal ideals in $\mathbf{R_p}$,   there is a smallest integer  $m$ such that $\mathbf{J} = \mathbf{I_1}^m ...\mathbf{I_s}^m$, and $\mathbf{R_p}= \bigoplus_j \mathbf{R_p}/\mathbf{I_j^m}$,  where each factor is a local ring.  $\square$

REMARK.  For finite, commutative, semisimple rings, several of the radical operators coalesce.
The radical mentioned in the theorem can be taken, for example, to be the intersection of maximal ideals, or as the nilpotent radical.

We use the term 'p splits' to mean that the core polynomial factors modulo(p), but that it does not ramify.

$\mathbf{THEOREM} \,6.5$

(6.51) If $p$ is inert,  then $\mathbf{R_p}$ is a field.

(6.52) If $p$ splits, then $\mathbf{R_p}$ has a trivial radical, thus is semisimple, i.e., is the direct sum of  $s$ simple rings (fields in this case), where $s$ is the number of prime ideals in the factorization of $\mathbf{R_p}$.

(6.53) If $p$ ramifies, then $\mathbf{R_p}$ has a non-trivial radical, and is a direct sum of  $s$ (non-trivial) local rings .  

PROOF.  Theorem 6.5 is a direct consequence of the structure theorem, Theorem 6.4. The $m$ in  theorem 6.4 is the l.c.m. of the ramification indices.   $\square$

REMARK. An ideal element $ \mathbf{m}$ outside of the radical is cyclic, i,e, satisfies $\mathbf{m}^n = \mathbf{m}$ for some natural number $n$. If $\mathbf{m} = \mathbf{e}$ is an idempotent, then the powers of $\mathbf{eA}$ coincide with the orbit of $\mathbf{e}$. This is because $(\mathbf{eA})^n = \mathbf{eA}^n$;  thus $\mathbf{eA}$ generates a cyclic group of order dividing $c_p[\mathbf{t}]$.
Using the standard matrix representation of elements in $\mathbf{R}$ or in $\mathbf{R_p}$, we can assign to each element a \textit{rank} by letting  $rank\mathbf{m} = rank\mathbf{M}$, where $\mathbf{M}$ is the standard matrix representation of $\mathbf{m}$.  We then observe that all elements in the same orbit have the same rank;  that the rank of a unit is $k$, the degree of the core polynomial; and,  that the rank of an ideal element is at 
 most the co-degree of the ideal , i.e., $k-d$, where $d$ is the degree of the minimal polynomial of the ideal. (The rank of the representing matrix cannot exceed the degree of the minimal polynomial). 

Denote the rank of an element $\mathbf{m}$ in $\mathbf{R_p}$ by $r(\mathbf{m})$.

$\mathbf{THEOREM \, 6.6}$

Suppose that $p$ splits and that $\{\mathbf{e_1},...,\mathbf{e_s}\}$ is a complete set of distinct primitive idempotents in $\mathbf{R_p}$. 
 $$r(\sum_1^s \mathbf{e_j}) =  \sum_1^s r(\mathbf{e_j}) = k.$$
 
Proof.  By (6.52),  $\mathbf{R_p}$ is semisimple.  We observe that:  $1 \le r(\mathbf{e_j}) < k$, and  since $\sum_1^s \mathbf{e_j} = 1$, $r(\sum_1^s \mathbf{e_j}) = k$.  The proof will then be a consequence of the following lemma and corollaries.

$ \mathbf{LEMMA}$ If we let $\mathbf{e}$ be the sum of the elements in any subset of the set of primitive idempotents $\{\mathbf{e_i}\}$, and let $\bar{\mathbf{e}}$ be the complementary sum, then $$r(\mathbf{e}) + r(\bar{\mathbf{e}})\leqslant k. $$ 
 
Proof.  If $\mathbf{E_1}$ and $\mathbf{E_2}$ are $ k \times k $ -matrices such that $\mathbf{E_1 E_2 }= \mathbf{0}$, then $r(\mathbf{E_1}) + r(\mathbf{E_2}) \leqslant k.$
 Since $\mathbf{E_1 E_2} = \mathbf{0}$ we have that $r(\mathbf{E_1}) \leqslant \nu(\mathbf{E_2})\leqslant k-r(\mathbf{E_2})$, where $\nu$ is the nullity of $\mathbf{E_2}$, and the lemma follows. $\square$

 $\mathbf{ COROLLARY \,6.61}$
\begin{center}
$r(\mathbf{e}) +  r(\mathbf{\bar{e}}) = k.$
\end{center}
Proof.  Using the above Lemma and the remark at the beginning of the proof of the theorem, we have $k = r(\mathbf{e} +  \mathbf{\bar{e}})  \leqslant   r(\mathbf{e}) +  r(\mathbf{\bar{e}}) = k.$  $\square$
\newpage  $\mathbf{ COROLLARY  \,6.62}$
 
 \begin{center}
 $r(\mathbf{e_i} + \mathbf{e_j}) = r(\mathbf{e_i}) + r(\mathbf{e_j}).$
\end{center}
 
 Proof.  From Corollary 6.61, we have that  $r(\mathbf{e_1}) +  r(\mathbf{\bar{e_1}}) = k$, so that we can apply the above arguments to  $r(\mathbf{\bar{e_1}} ) = k -  r(\mathbf{e_1}) $ to deduce that $r(\sum_2^s\mathbf{e_i} ) = \sum_2^s r(\mathbf{e_i})$; hence, Corollary 6.62 holds.  $\square$

 Theorem 6.6 follows now from the proof of Corollary 6.62. $\square$

$\mathbf{COROLLARY}$ 6.63  

\noindent   If we let $\mathbf{B_1}, ..., \mathbf{B_s}$ be the ideals  $\mathbf{R_p }\mathbf{\mathbf{e_1}}, ..., \mathbf{R_p} \mathbf{e_s}$ in $\mathbf{R_p}$, and let $\mathbf{B_j}^* = \mathbf{B_j }- \{\mathbf{0}\}$, then 
 $$\mathbf{B_1}^*\times...\times \mathbf{B_s}^* = \mathbf{\mathbf{G_p}},$$
 where $\mathbf{G_p}$ is the group of units of $\mathbf{R_p}$.
 If  $p$ does not ramify, the $\mathbf{B_j}^*$ are finite fields.

Proof.  This follows from Theorem 6.5, Theorem 6.6, the fact that ranks of non-zero elements of $\mathbf{R_p}$ are positive integers,  and that an element of $\mathbf{R_p}$ is a unit if and only if its norm is not zero  \cite{McD}. $\square$

$\mathbf{COROLLARY \,6.64}$
$$|\mathbf{G_p}| = |\mathbf{B_1}^*| ... |\mathbf{B_s}^*| = (p^{r_1}-1)...(p^{r_s}-1),$$ where $p^{r_i}$ is the order of $\mathbf{B_i}$ and $r_i$ is the rank of $\mathbf{e_i}$.     $\square$
  
$\mathbf{COROLLARY \,6.65}$

If $p$ splits,  the period $c_p[t_1,...,t_k] = lcm \{c_p(\mathbf{minpoly(e_i}))\}_1^s $. $\square$

The following result gives a remarkable connection between the p-periodicity of a linear recursion and the splitting properties of primes in associated rational number fields.
\newpage $\mathbf{THEOREM\,6.7}$   $p$ divides $c_p[\mathbf{t}]$ if and only if $p$ ramifies.

Proof.  First, we assume that $p | c_p[\mathbf{t}]$ and that $p$ does not ramify;  but then, by Theorem 6.5, $\mathbf{R_p}$ is semisimple, and, so by Corollary 6.64, $p | (p^{r_i}-1)$ for some $i$. A contradiction.
In particular, If $p$ does not ramify, $|\mathbf{G_p}|$ and $p$ are relatively prime.
In order to prove the converse,  we first prove the following lemma:

$\mathbf{LEMMA} $  If $\mathbf{e}$ is an idempotent in an ideal of $\mathbf{R_p}$, then the $\mathbf{H_p}$-orbit of $\mathbf{e}$ consists of the powers of $\mathbf{eA}$, a multiplicative cyclic group. In particular, the order of $\mathbf{eA}$ divides $c_p[\mathbf{t}]$.

Proof.  All of this follows easily from the fact that $ (\mathbf{e A)^{n}}   =  \mathbf{e}^{n} \mathbf{A^{n}}  = \mathbf{e A^n}$, that $\mathbf{eA^{c_p}} =\mathbf{ e}$,  and that the length of any orbit divides the period. $\square $

To finish the proof of the theorem, we observe that,  if $p$ ramifies, then each of the direct factors in $\mathbf{R_p}$ is a non-trivial local ring, say $\mathbf{B_j}$, where $\mathbf{B_j} =  \mathbf{I_j}/ {\mathbf{I_j}}^m$. If $\mathbf{B_j}^*$ is the group of units in $\mathbf{B_j}$,  then there is an idempotent  $\mathbf{e_j}$ in $\mathbf{I_j}$ and  $\mathbf{e_j }+ \mathbf{m}$ is a unit in the local ring $\mathbf{B_j}$, that is, is in $\mathbf{B_j}^{*}$,   whenever $\mathbf{m }\in {\mathbf{I}_\mathbf{j}^m}$, i.e. whenever $\mathbf{m }\in \mathbf{J}$.
Moreover, there is a bijective correspondence between such $\mathbf{m}'s$ in $\mathbf{I_j}$ and the elements in the orbit of $\mathbf{e_j}$, so that  by the Lemma, $p$ divides $c_p$ .  Thus $p$ divides $|\mathbf{H_p}|$  and, hence, the order of $\mathbf{G_p}$. $\square$

 \textbf{REMARK}: The well-known fact that a rational prime $p$ ramifies with respect to a cyclotomic extension over $CP(n)$  only if $p$ divides  $n$ now follows immediately from Theorems 2.1, 4.14, 6.7.
 
\textbf{ REMARK}: Note that our notation $c_p[\textbf{t}]$ for period of the recursion determined by the core polynomial $[t_1,...,t_k]$,  for $p=1$ or $p$ prime, could as well be written  $c_p(\mathbf{R_p})$,  where, of course,  
 $\mathbf{R_p} = \mathbf{R}$  if $p=1$.  (Here it really doesn't matter whether the core is irreducible or not, that is,  whether $\mathbf{R}$ is a number field,  or merely a commutative ring.) Since each maximal ideal in $\mathbf{R_p}$ is determined by and determines its minimal polynomial, which in turn determines the periods of the p-factors of the core polynomial of  $\mathbf{R}$,  the notation $c_p(\mathbf{I})$ can be used to denote the period of the recursion determined by a factor of the original core polynomial.  We use this notation in the statement of the next theorem. The structure theory for semilocal rings
now provides the right setting in which to reformulate and prove Proposition 5.1.  We do this in 

$\mathbf{THEOREM \,6.8}$

Suppose that  the semilocal ring $\mathbf{R_p}$ has maximal ideals $\mathbf{I_1,...,I_s}$, that is,  suppose that the core polynomial has $s$ irreducible factors,  denoting the radical of $\mathbf{R_p}$ by $ \mathbf{J}$, then,

(1)  $$|\mathbf{G_p(R_p})| = |\mathbf{B_1^*}| \cdot...\cdot |\mathbf{B_s^*}|\cdot|\mathbf{J}|,$$

(2)  $$c_p(\mathbf{R_p}) = lcm\{c_p(\mathbf{I_1}),...,c_p(\mathbf{I_s})\} \cdot |\mathbf{J}| ,$$ 

where $\mathbf{B_j^*} = \mathbf{I_j}/\mathbf{J},  j=1,...,s,$ and $p=1$ or $p$ is prime.

Proof.  Consider the exact sequence  $$ \mathbf{J }\rightarrowtail \mathbf{R_p}  \twoheadrightarrow \frac{\mathbf{R_p}}{\mathbf{J}}  $$ where  $ \mathbf{J}$ is the radical of the ring. Since $\frac{\mathbf{R_p}}{\mathbf{J}}$ is semi-simple, by Corollary 6.65, $c_p(\frac{\mathbf{R_p}}{\mathbf{J}})= lcm \{c_p(\mathbf{minpoly(\mathbf{B_j^*})})\}_1^s $. If $\mathbf{R_p}$ splits, that is, if the radical is  $\mathbf{0}$, then we are done.  In any case,  by Corollary 6.64,  $|\mathbf{G_p}(\mathbf{R_p/J})| =  |\mathbf{B_1^*}| \times...\times|\mathbf{B_s^*}|\times|\mathbf{J}|,.$ It is clear that units in $\mathbf{R_p}$ are mapped homomorphically onto the units of $\mathbf{R_p/J},$  and since  $\mathbf{u} \rightarrow \mathbf{u+J}, \mathbf{u} \in \mathbf{G_p}$, part (1) of the theorem follows.

As for part (2) of the theorem,  we have by Corollary 6.65 that $c_p(\mathbf{R_p/J}) = lcm\{ |\mathbf{\mathbf{B_1^*}}|,..., |\mathbf{B_s^*}|\}$ = $lcm\{(p_1^{r_1}-1),...,(p_1^{r_s}-1)\}$ = $lcm \{c_p(\mathbf{I_j})\}$, that is, the least common multiple of the periods of the irreducible factors of the core polynomial of $\mathbf{R_p/J}$.  But this number is just the order of the period group $(\mathbf{H_p + J})/\mathbf{J}$ in  $\mathbf{R/J}$, which accounts for the factor $|\mathbf{J}|$ in part (2) of the theorem.
 $\square$

 Theorem 6.8 now enables us to prove a detailed version of Proposition 5.1:
\newpage  
$\textbf{ COROLLARY 6.9}$

Suppose the core polynomial (reducible or not) factors mod(p) into $s$ irreducible factors, then the p-core polynomial is the least common multiple of the periods of the irreducible factors times the order of the radical of the semi-simple ring $\mathbf{R_p}$.  Furthermore,  $\mathbf{J}$ is non-trivial exactly when  $p$  divides the period.  In terms of algebraic number fields,  this is just the case when  $p$ ramifies, i.e. when $p$ divides the discriminant of the field.   $\square$

The following example illustrates the situation in the case of a non-trivial radical:

\begin{figure}[h]
\centering
\includegraphics[width=3in]{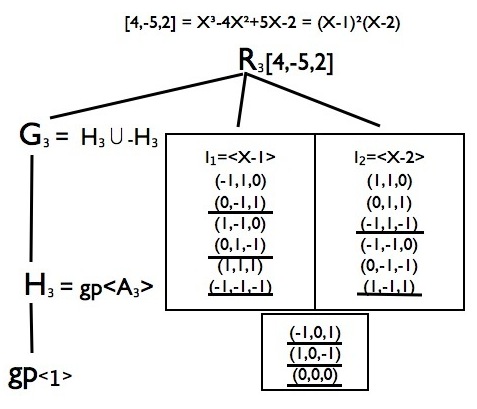} 
\caption{{\bf non-trivial radical}}
\label{Figure 3.}
\end{figure}
 
Note that $(0,-1,-1)$ is a unit in $\mathbf{I_2}$ and that $\mathbf{I_2-J}$ is a multiplicative group in which this unit acts as the identity.  In this example  $|\mathbf{G_3}| = 12$ and the \textit{Period subgroup} = $gp<\mathbf{A}> = \mathbf{H_3}$ has order 6,
$|\mathbf{H_3}| = 6$

In the case of number fields, it is a trivial fact that only finitely many primes ramify. There are examples however when,  for every  $p$, $p$  divides the period;  for instance, the following is such a case.  Consider the core polynomial $[2,-1]$.  The 2-linear sequence $\{F_n[2,-1]\}$ for this core is just
$\{1,2,3,...,n,...\}$.  It is easy to see that $c(\mathbf{R_p}) = p$ for every  $p$.

  Key Words:  Symmetric polynomials, linear recursive sequences, number fields.
 
  Trueman MacHenry \\
  York University, Toronto, Canada \\
  machenry@mathstat.yorku.ca \\
  
  Kieh Wong \\
  Centennial College, Toronto, Canada \\
  kkwong@centennialcollege.ca
  
   May 9, 2007
  \end{document}